# INVERSE SCATTERİNG PROBLEM ON THE HALF-AXİS FOR A FİRST ORDER SYSTEM OF ORDİNARY DİFFERENTİAL EQUATİONS

MANSUR I. ISMAİLOV

Gebze Institute of Technology, Turkey, Gebze-Kocaeli
E mail: mismailov@gyte.edu.tr

**ABSTRACT.** In this article, the inverse scattering problem (ISP) of recovering the matrix coefficient of a first order system of ordinary differential equations on the half-axis from its scattering matrix is considered. In the case of a triangular structure of the matrix coefficient, this system has a Volterra-type integral transformation operator at infinity. Such type of transformation operator allows to determine the scattering matrix on the half-axis via the matrix Riemann-Hilbert factorization in the case, where contour is real axis, normalization is canonical and all the partial indices are zero. The ISP on the half-axis is solved by reducing it to ISP on the whole axis for the considered system with the coefficients that are extended to the whole axis as zero.

## 1. INTRODUCTİON

The inverse scattering problems (ISP's) for the differential equations consists in reconstructing these equations from their scattering operators. Such problems appear in mathematics, mechanics and various areas of natural sciences and engineering. Various settings of ISP's for systems of ordinary differential equations were examined by many authors. Many publications are devoted to the ISP for the Dirac systems and their $2n$ generalizations (see [1–7] and the references therein). The main approach in this field is based on the transformation operator (TO) method, and the results proved in this way are analogs of results obtained for the Sturm–Liouville operator [8–10].

We consider the following system of ODE on the half-axis:

$$l(y) = -i\frac{dy}{dx} + Q(x)y = \lambda \sigma y, \ 0 \leq x < +\infty \tag{1.1}$$

with the parameter $\lambda$. It is assumed that $\sigma = \begin{bmatrix} \sigma_1 & 0 \\ 0 & \sigma_2 \end{bmatrix}$ is $2n \times 2n$ order diagonal matrix with constant diagonal elements, where $\sigma_1 = diag(\xi_1, \ldots, \xi_n)$, $\sigma_2 = diag(\xi_{n+1}, \ldots, \xi_{2n})$, $\xi_1 < \ldots < \xi_n < 0 < \xi_{n+1} < \ldots < \xi_{2n}$, and $Q$ is $2n \times 2n$ order matrix function with measurable complex valued entries, which the matrix norm satisfies the condition

$$\|Q(x)\| \leq Ce^{-\varepsilon x} \tag{1.2}$$

where $C, \varepsilon$ are positive constants. The matrix $Q(x)$ is called as potential.

Let $y(x,\lambda) = \begin{bmatrix} y_1(x,\lambda) \\ y_2(x,\lambda) \end{bmatrix}$, where $y_1(x,\lambda)$ and $y_2(x,\lambda)$ are $n$ dimensional vector functions. Consider the system (1.1) under the boundary condition at point $x = 0$ in following form:

$$y_2(0,\lambda) = H y_1(0,\lambda), \quad det H \neq 0. \tag{1.3}$$

The system (1.1) is important from an application point of view. It is used in solving N-Wave problem that finds application in nonlinear optics ([11], p. 174).

The ISP on the whole axis for the operator (1.1) has been treated in many papers. We recall the important papers [12- 14] (see also [11], p. 163 and [15], p. 7), where the ISP for system (1.1) on the whole axis was investigated by using the Riemann - Hilbert problem for analytic functions. The ISP for system (1.1) on the half-axis in the cases $n > 1$, have not been satisfactorily investigated yet. The principal difficulty is to determine the sufficient quantity of scattering problems (on the semi-axis for the same system) ensuring the uniqueness of the solution of the ISP under consideration (see Example in present paper). The ISP on the half-axis is closely connected with the ISP on the whole axis. Such type a relationship is demonstrated in [16], where the ISP for the first order hyperbolic system on the half-axis is investigated.

In contrast to the inverse scattering problem, the inverse problems problems for the system (1.1) from spectral function and from Weyl function were satisfactorily investigated. We recall the papers [17] and [18], where the inverse problem for the system (1.1) on the half-axis and on the finite interval with self-adjoint potential was investigated from its spectral function by using TO method, and the papers [19, 20] where the inverse problem for the system (1.1) on the half-axis was studied from its matrix Weyl function.

In this paper, we consider the system (1.1) with the matrix potential which has some triangular structure. We show that this system has Volterra type integral TO at infinity. The integral representation of the solution of the system (1.1), which is taken from TO at infinity, allows to determine the scattering matrix on the half-axis and some its analytical properties. From the uniqueness point of view, we consider ISP from two scattering matrix on the semi-axis. Under the condition of unique solvability of some matrix Riemann-Hilbert problem with zero index, the ISP on the half-axis is reduced to ISP on the whole axis for the same system with the coefficients which are zero for $x < 0$. However, for the general system (1.1) with arbitrary potential $Q(x)$, the ISP turns out to be substantially more difficult. In the general case, the transformation operators have a rather complicated structure and this makes the solution of the inverse problem difficult.

The paper is organized as follows: In next section we determine Volterra type integral TO for the system (1.1) at infinity, when the matrix coefficients $Q(x)$ of the system are in the special triangular structure. In Section 3 the scattering matrix on the half-axis is defined and some its analytic properties are studied. In Section 4 the ISP for the system (1.1) on the half-axis is formulated and the uniqueness of the solution of the problem is proved by reducing that problem to the ISP on the whole axis for the system with the coefficients which are zero for $x < 0$. In Section 5 the non-uniqueness in the ISP is discussed, when there are some violations on the conditions. Finally, in Appendix some conditions for the unique solvability of the matrix Riemann-Hilbert (RH) problem when the contour is real axis, normalization is canonical and all the partial indices are zero, are given.

**DEFINITION.** We will call the system (1.1) as an $\mathcal{M}$ -canonical system, if the matrix potential $Q = \begin{bmatrix} q_{11} & q_{12} \\ q_{21} & q_{22} \end{bmatrix}$ of the system (1.1) are given by the following triangular structure:

- $q_{11}$ is a strictly lower triangular $n \times n$ matrix function;
- $q_{12}$ is a lower anti-triangular $n \times n$ matrix function;
- $q_{21}$ is an upper anti-triangular $n \times n$ matrix function;
- $q_{22}$ is a strictly upper triangular $n \times n$ matrix function.

The system (1.1) can be consider as the stationary analogue of the first order system of hyperbolic differential equations. As regards different inverse problems for first order system of partial differential equations are studied in monographs [20, 21] and the literature cited in there. The uniqueness of the solution of the ISP for the nonstationary $\mathcal{M}$-canonical system (1.1) was studied in [23].

**NOTATION.** Throughout the paper, we shall write $\mathbb{R}$ for the real axis, $I$ for $n \times n$ order identity matrix. The upper and lower complex $\lambda$-plane denote the sets with $Im\lambda \geq 0$ and $Im\lambda \leq 0$, respectively. The matrix $(a_{ij})_{i,j=1}^{n}$ is called lower (strictly) triangular, if $a_{ij} = 0$ for all $j > i$ (for all $j \geq i$); lower anti-triangular, if $a_{ij} = 0$ for all $j + i \leq n$; upper (strictly) triangular, if $a_{ij} = 0$ for all $j < i$ (for all $j \leq i$); and upper anti-triangular, if $a_{ij} = 0$ for all $j + i \geq n + 2$.

## 2. THE TRANSFORMATİON OPERATOR FOR THE SYSTEM (1.1) WİTH THE BOUNDARY CONDİTİON AT İNFİNİTY

In solving ISPs the Volterra type integral representation of the solution plays an important role. Such a representation for the $\mathcal{M}$-canonical system (1.1) can be taken from a transformation operator for that system with the following boundary condition at infinity

$$y_1(x, \lambda) = e^{i\lambda\sigma_1 x} A + o(1), \quad x \to +\infty,$$

$$y_2(x, \lambda) = e^{i\lambda\sigma_2 x} B + o(1), \quad x \to +\infty, \quad (2.1)$$

where $e^{i\lambda\sigma_1 x} = diag(e^{i\lambda\xi_1 x}, \ldots, e^{i\lambda\xi_n x})$, $e^{i\lambda\sigma_2 x} = diag(e^{i\lambda\xi_{n+1} x}, \ldots, e^{i\lambda\xi_{2n} x})$.

More precisely, if the system (1.1) is $\mathcal{M}$-canonical then it has the solution in form of Volterra integral operator at infinity.

We first prove the following lemma.

**LEMMA 1.** *Let $\lambda$ be a real number and the potential $Q(x)$ satisfy the condition (1.2), then the following statements hold for the system (1.1).*

1) *There exists a unique bounded solution of the system (1.1) with the boundary condition (2.1);*
2) *For a bounded solution $y(x, \lambda) = \begin{bmatrix} y_1(x, \lambda) \\ y_2(x, \lambda) \end{bmatrix}$ of the system (1.1), the asymptotic relations (2.1) are true, where*

$$A = y_1(0, \lambda) - i \int_0^{+\infty} e^{-i\lambda\sigma_1 s} [q_{11}(s)y_1(s, \lambda) + q_{12}(s)y_2(s, \lambda)]ds,$$

$$B = y_2(0,\lambda) - i\int_0^{+\infty} e^{-i\lambda\sigma_2 s}[q_{21}(s)y_1(s,\lambda) + q_{22}(s)y_2(s,\lambda)]ds. \qquad (2.2)$$

**PROOF.** Let $\lambda$ be real number. The first statement of the lemma is equivalent to unique solvability of the following system of integral equation in class of bounded function:

$$y_1(x,\lambda) = e^{i\lambda\sigma_1 x}A + i\int_x^{+\infty} e^{i\lambda\sigma_1(x-s)}[q_{11}(s)y_1(s,\lambda) + q_{12}(s)y_2(s,\lambda)]ds,$$

$$(2.3)$$

$$y_2(x,\lambda) = e^{i\lambda\sigma_2 x}B + i\int_x^{+\infty} e^{i\lambda\sigma_2(x-s)}[q_{21}(s)y_1(s,\lambda) + q_{22}(s)y_2(s,\lambda)]ds.$$

Applying the method of successive approximation to the system (2.3) and using the condition (1.2) we obtain that the solution of this system exists and is unique in class of bounded function.

Let us prove the second statement of the lemma.

Let $y(x,\lambda) = \begin{bmatrix} y_1(x,\lambda) \\ y_2(x,\lambda) \end{bmatrix}$ is a bounded solution of the system (1.1). Since the system (1.1) have the fundamental system of bounded solution in the case of real $\lambda$ (see Theorem 8.1. in [24]), then we have

$$y_1(x,\lambda) = e^{i\lambda\sigma_1 x}y_1(0,\lambda) - i\int_0^x e^{i\lambda\sigma_1(x-s)}[q_{11}(s)y_1(s,\lambda) + q_{12}(s)y_2(s,\lambda)]ds,$$

$$y^2(x,\lambda) = e^{i\lambda\sigma_2 x}y^2(0,\lambda) - i\int_0^x e^{i\lambda\sigma_2(x-s)}[q_{21}(s)y^1(s,\lambda) + q_{22}(s)y^2(s,\lambda)]ds,$$

by the uniqueness of the solution of Cauchy problem for the system (1.1) at the point $x = 0$. Hence, when $\lambda$ is real, the solution $y(x,\lambda)$ of the system (1.1) has the asymptotic form (2.1) as $x \to +\infty$, where $A$ and $B$ are determined as in (2.2). Lemma 1 is proved.

Now, suppose that the system (1.1) has the $\mathcal{M}$-canonical form. This system with the boundary condition (2.1) has the solution in the form

$$y_1(x,\lambda) = e^{i\lambda\sigma_1 x}A + \int_x^{+\infty} A_{11}(x,t)e^{i\lambda\sigma_1 t}Adt + \int_x^{+\infty} A_{12}(x,t)e^{i\lambda\sigma_2 t}Bdt,$$

$$(2.4)$$

$$y_2(x,\lambda) = e^{i\lambda\sigma_2 x}B + \int_x^{+\infty} A_{21}(x,t)e^{i\lambda\sigma_1 t}Adt + \int_x^{+\infty} A_{22}(x,t)e^{i\lambda\sigma_2 t}Bdt,$$

where $A_{11}(x,t)$ is a lower triangular $n \times n$ matrix, $A_{12}(x,t)$ is a lower anti-triangular $n \times n$ matrix, $A_{21}(x,t)$ is an upper anti-triangular $n \times n$ matrix and $A_{22}(x,t)$ is an upper triangular $n \times n$ matrix, by starting (2.3) and (2.4) we obtain the system of integral equation with respect to matrix kernel $A_{ij}(x,t)$ $(i,j = 1,2)$:

$$[A_{11}(x,t)]_{kj} = i\frac{\xi_j}{\xi_j-\xi_k}[q_{11}]_{kj}(x+\frac{\xi_j}{\xi_j-\xi_k}(t-x))$$

$$+i\int_x^{x+\frac{\xi_j}{\xi_j-\xi_k}(t-x)}[q_{11}A_{11}+q_{12}A_{21}]_{kj}(s,t+\frac{\xi_k}{\xi_j}(s-x))ds,$$

(2.5)

$$[A_{21}(x,t)]_{kj} = i\frac{\xi_j}{\xi_j-\xi_{n+k}}[q_{21}]_{kj}(x+\frac{\xi_j}{\xi_j-\xi_{n+k}}(t-x))$$

$$+i\int_x^{x+\frac{\xi_j}{\xi_j-\xi_{n+k}}(t-x)}[q_{21}A_{11}+q_{22}A_{21}]_{kj}\left(s,t+\frac{\xi_{n+k}}{\xi_j}(s-x)\right)ds.$$

$$t \geq x.$$

$$[A_{12}(x,t)]_{kj} = i\frac{\xi_{n+j}}{\xi_{n+j}-\xi_k}[q_{12}]_{kj}(x+\frac{\xi_{n+j}}{\xi_{n+j}-\xi_k}(t-x))$$

$$+i\int_x^{x+\frac{\xi_{n+j}}{\xi_{n+j}-\xi_k}(t-x)}[q_{12}A_{22}+q_{11}A_{12}]_{kj}(s,t+\frac{\xi_k}{\xi_{n+j}}(s-x))ds,$$

(2.6)

$$[A_{22}(x,t)]_{kj} = i\frac{\xi_{n+j}}{\xi_{n+j}-\xi_{n+k}}[q_{22}]_{kj}(x+\frac{\xi_{n+j}}{\xi_{n+j}-\xi_{n+k}}(t-x))$$

$$+i\int_x^{x+\frac{\xi_{n+j}}{\xi_{n+j}-\xi_{n+k}}(t-x)}[q_{21}A_{12}+q_{22}A_{22}]_{kj}(s,t+\frac{\xi_{n+k}}{\xi_{n+j}}(s-x))ds,$$

$$t \geq x,$$

where $k,j = 1,\ldots,n$; $[A]_{kj}$ denotes the $k,j$ element of the matrix function $A$; $[qA]_{kj}(s,t) \equiv [q(s)A(s,t)]_{kj}$; when $k=j$ the upper bound of the integral is $+\infty$.

Applying the method of successive approximation to the systems (2.5) and (2.6) we obtain that the solution of these systems exist and are unique in class of bounded function, here $A_{ij}(x,t)$, $t \geq x \geq 0$, $i,j = 1,2$, has the estimate

$$\|A_{ij}(x,t)\| \leq \tilde{C}e^{-\varepsilon(x+\theta(t-x))},$$

(2.7)

where $\tilde{C}$ is constant, $\theta = min\{\theta_1,\theta_2,\theta_3,\theta_4\}$, $\theta_1 = min_{k>j}\frac{\xi_j}{\xi_j-\xi_k}$, $\theta_2 = min_{k+j>n}\frac{\xi_{n+j}}{\xi_{n+j}-\xi_k}$, $\theta_3 = min_{k+j<n+2}\frac{\xi_j}{\xi_j-\xi_{n+k}}$ и $\theta_4 = min_{k<j}\frac{\xi_{n+j}}{\xi_{n+j}-\xi_{n+k}}$.

When we let $t = x$ we obtain from (2.4) and (2.5) that

$$[q_{11}(x)]_{kj} = -i\frac{\xi_j-\xi_k}{\xi_j}[A_{11}(x,x)]_{kj}, \qquad [q_{12}(x)]_{kj} = -i\frac{\xi_{n+j}-\xi_k}{\xi_{n+j}}[A_{12}(x,x)]_{kj},$$

$$[q_{21}(x)]_{kj} = -i\frac{\xi_j-\xi_{n+k}}{\xi_j}[A_{21}(x,x)]_{kj}, \qquad [q_{22}(x)]_{kj} = -i\frac{\xi_{n+j}-\xi_{n+k}}{\xi_{n+j}}[A_{22}(x,x)]_{kj}.$$

(2.8)

Now, we prove the following lemma.

**LEMMA 2.** *Let $\lambda$ be a real number and $y(x,\lambda) = \begin{bmatrix} y_1(x,\lambda) \\ y_2(x,\lambda) \end{bmatrix}$ be a bounded solution of the $\mathcal{M}$ - canonical system (1.1) with the potential $Q(x)$ satisfying the condition (1.2). Then the representation (2.4) of the solution holds, where the matrix kernels $A_{ij}(x,t)$ $(i,j = 1,2)$ satisfy the system of integral equations (2.5) and (2.6). In addition, the estimation (2.7) and equality (2.8) hold.*

**PROOF.** Let $y(x,\lambda)$ be a bounded solution of $\mathcal{M}$ - canonical system (1.1). Then the asymptotic relation (2.1) holds according to the second statement of Lemma 1. It is clear that the solution $y(x,\lambda)$ of the system (1.1) with the boundary condition (2.1) satisfies the system of integral equations (2.3). When the function $Q(x)$ satisfies the estimate (1.2) then there exists a unique solution of the system of integral equations (2.5) and (2.6) in the class of bounded functions, since the mentioned systems are Volterra type integral equations. In addition, the estimation (2.7) in $t \geq x \geq 0$ and the equalities (2.8) are valid for these solutions. Conversely, if the functions $A_{ij}(x,t)$ $(i,j = 1,2)$ satisfy the systems (2.5) an (2.6), then the formula (2.4) gives the bounded solution of the system (2.3) for real $\lambda$. Since the bounded function $y(x,\lambda)$ satisfies the system (2.3) and the system (2.3) has a unique solution in the class of bounded functions, $y(x,\lambda)$ is represented in the form of (2.4). Lemma 2 is proved.

## 3. SCATTERİNG MATRİX AND İTS PROPERTİES

Suppose that, for real $\lambda$, $y(x,\lambda) = \begin{bmatrix} y_1(x,\lambda) \\ y_2(x,\lambda) \end{bmatrix}$ is a bounded solution of the $\mathcal{M}$ -canonical. Then by Lemma 2 the solution $y(x,\lambda)$ is representable in the form (2.4). Hence from (2.4) and (1.3) we have

$$A_{21-}(\lambda)A + (I + A_{22+}(\lambda))B = H[(I + A_{11-}(\lambda))A + A_{12+}(\lambda)B], \tag{3.1}$$

where

$$A_{k1-}(\lambda) = \int_0^{+\infty} A_{k1}(0,t)e^{i\lambda\sigma_1 t}\,dt, \qquad A_{k2+}(\lambda) = \int_0^{+\infty} A_{k2}(0,t)e^{i\lambda\sigma_2 t}\,dt.$$

By denoting

$$A_{H+}(\lambda) = A_{22+}(\lambda) - HA_{12+}(\lambda),$$
$$A_{H-}(\lambda) = HA_{11-}(\lambda)H^{-1} - A_{21-}(\lambda)H^{-1}.$$

the formula (3.1) has the form

$$[I + A_{H+}(\lambda)]B = [I + A_{H-}(\lambda)]HA. \tag{3.2}$$

We introduce the matrix function

$$S_H(\lambda) = [I + A_{H+}(\lambda)]^{-1}[I + A_{H-}(\lambda)], \quad \lambda \in \mathbb{R}. \tag{3.3}$$

By analogy with the case $n = 1$ (see [4]) we call $S_H(\lambda)$, $\lambda \in \mathbb{R}$ the scattering matrix on the half-axis for the $\mathcal{M}$-canonical system with the boundary condition (1.3).

The following lemma is true.

**LEMMA 3.** *The matrix functions $I + A_{H-}(\lambda)$ u $I + A_{H+}(\lambda)$ are analyitic for $Im\lambda < -\frac{\theta}{\xi_1}\varepsilon$ and $Im\lambda > -\frac{\theta}{\xi_{2n}}\varepsilon$, respectively. The following asymptotic relations also hold as $|\lambda| \to +\infty$:*

$$det[I + A_{H-}(\lambda)] = 1 + o(1), \quad \left(Im\lambda < -\frac{\theta}{\xi_1}\varepsilon\right),$$

$$det[I + A_{H+}(\lambda)] = 1 + o(1), \quad \left(Im\lambda > -\frac{\theta}{\xi_{2n}}\varepsilon\right). \tag{3.4}$$

**PROOF.** From the estimation (2.7) of the matrix kernels $A_{ij}(x,t)$ $(i,j = 1,2)$ follows that $\|A_{ij}(0,t)\| \leq \tilde{C}e^{-\varepsilon\theta t}$. It means that the matrix functions $A_{k1-}(\lambda) = \int_0^{+\infty} A_{k1}(0,t)e^{i\lambda\sigma_1 t} dt$ and $A_{k2+}(\lambda) = \int_0^{+\infty} A_{k2}(0,t)e^{i\lambda\sigma_2 t} dt$ are analytic for $Im\lambda < -\frac{\theta}{\xi_1}\varepsilon$ and $Im\lambda > -\frac{\theta}{\xi_{2n}}\varepsilon$, respectively, and tend to zero as $|\lambda| \to +\infty$ in the domains of analyticity. Lemma 3 is proved.

According to Lemma 3, the functions $det[I + A_{H+}(\lambda)]$ and $det[I + A_{H-}(\lambda)]$ have a finite number of zeros. We define $\varepsilon_0 > 0$ by the relation

$$\varepsilon_0 = \min\left\{\varepsilon_1, -\frac{\theta}{\xi_1}\varepsilon, \frac{\theta}{\xi_{2n}}\varepsilon\right\},$$

where $\varepsilon_1$ be the distance from the real axis to the non-real zeros of the functions $det[I + A_{H+}(\lambda)]$ and $det[I + A_{H-}(\lambda)]$. Then the relations

$$det[I + A_{H+}(\lambda)] \neq 0, \quad det[I + A_{H-}(\lambda)] \neq 0 \tag{3.5}$$

hold for $0 < |Im\lambda| < \varepsilon_0$.

The next theorem about the properties of scattering matrix $S_H(\lambda)$ follows from (3.5) and Lemma 3.

**THEOREM 1.** *The matrix functions $S_H(\lambda)$ and $S^{-1}{}_H(\lambda)$ are meromorphic in the strip $|Im\lambda| < \varepsilon_0$, they have not nonreal poles and as $|\lambda| \to +\infty$*

$$S_H(\lambda) = I + o(1), \quad S^{-1}{}_H(\lambda) = I + o(1).$$

### 4. ISP ON THE HALF-AXİS AND MATRİX RİEMANN-HİLBERT PROBLEMS CORRESPONDİNG TO İT

The inverse scattering problem (ISP) on the half-axis for the system (1.1) consists in recovering the matrix potential $Q(x)$ from a given matrix function $S_H(\lambda)$. The exact solvable examples show that one scattering problem is not enough for the unique restoration of the potential. In this way, we consider the restoration of the potential from two scattering matrices which correspond to different boundary conditions in form of (1.3).

Let $S_{H_1}(\lambda)$ and $S_{H_2}(\lambda)$ be two scattering matrices on the half-axis for the $\mathcal{M}$-canonical system (1.1), where
$$\det(H_1 - H_2) \neq 0. \tag{4.1}$$

By the definition of scattering matrix on the half-axis we get
$$[I + A_{H_k+}(\lambda)]S_{H_k}(\lambda) = [I + A_{H_k-}(\lambda)], \quad \lambda \in \mathbb{R} \quad k = 1,2, \tag{4.2}$$

where
$$A_{H_k+}(\lambda) = A_{22+}(\lambda) - H_k A_{12+}(\lambda),$$
$$A_{H_k-}(\lambda) = H_k A_{11-}(\lambda) H_k^{-1} - A_{21-}(\lambda) H_k^{-1}. \tag{4.3}$$

If these matrices are known, then relations (4.2) are matrix Riemann-Hilbert problems, where the contour is real axis and normalization is canonical. We will call these problems as Riemann-Hilbert problems of the ISP on the half-axis for the $\mathcal{M}$-canonical system (1.1).

The ISP on the half-axis for the system (1.1) closely is related with the ISP on the whole axis. In this reason, we introduce the matrix $P(\lambda)$, $\lambda \in$ for the bounded solutions $y(x,\lambda)$ of the M-canonical system as follows
$$P(\lambda)\begin{bmatrix}A\\B\end{bmatrix} = \begin{bmatrix}y_1(0,\lambda)\\y_2(0,\lambda)\end{bmatrix}, \quad \lambda \in \mathbb{R}. \tag{4.4}$$

By the uniqueness of the solution of the Cauchy problem at the point $x = 0$, for the system (1.1) we have $y(x,\lambda) = 0$ when $y_1(0,\lambda) = y_2(0,\lambda) = 0$. Then $A = B = 0$, by the formula (2.2). It means that the matrix $P(\lambda)$ is invertable. We will call the matrix $\Pi(\lambda) = P^{-1}(\lambda)$ the transmission matrix.

Now, consider the system of ODE on the whole axis
$$-i\frac{dy}{dx} + \tilde{Q}(x)y = \lambda \sigma y, \tag{4.5}$$

with the potential $\tilde{Q}(x) = \begin{cases} Q(x), & x \geq 0 \\ 0, & x < 0 \end{cases}$.

By comparing the definition (4.4) of the transmission matrix $\Pi(\lambda)$ with the definition of scattering matrix on the whole axis (see [12]), it is clear to see that matrix $\Pi(\lambda)$ is the scattering matrix for the system (4.5). The ISP for such type of systems on the whole axis is studied in [12-15] (see also [11], p. 163). Thus, we obtain the following result about the ISP for the $\mathcal{M}$-canonical system (1.1) on the half-axis.

**THEOREM 2.** *Let $S_{H_1}(\lambda)$ and $S_{H_2}(\lambda)$ be two scattering matrices on the half-axis for the $\mathcal{M}$-canonical system (1.1) with the potential $Q(x)$ satisfying the condition (1.2). Let the matrices $H_1$ and $H_2$ satisfy the condition (4.1). Then, the matrix $Q(x)$ is uniquely determined from matrices $S_{H_1}(\lambda)$ and $S_{H_2}(\lambda)$, when the Riemann-Hilbert problems (4.2) are uniquely solvable.*

**PROOF.** Once, let us show that the transmission matrix $\Pi(\lambda)$ is uniquely determined from $S_{H_1}(\lambda)$ and $S_{H_2}(\lambda)$, when the Riemann-Hilbert problems (4.2) are uniquely solvable. Applying the representation (2.4) to definition (4.4) of $P(\lambda)$ it can be easily seen the following block structure of $P(\lambda)$:

$$P(\lambda) = \begin{bmatrix} I+A_{11-}(\lambda) & A_{12+}(\lambda) \\ A_{21-}(\lambda) & I+A_{22+}(\lambda) \end{bmatrix}.$$

When the Riemann-Hilbert problems (4.2) are uniquely solvable, i.e. if matrices $A_{H_k+}(\lambda)$ and $A_{H_k-}(\lambda)$ are uniquely determined from $S_{H_k}(\lambda)$ ($k=1,2$) in domain which they are analytical, then matrices $A_{11-}(\lambda)$, $A_{21-}(\lambda)$, $A_{12+}(\lambda)$ and $A_{22+}(\lambda)$ are uniquely expressed from $A_{H_k+}(\lambda)$ and $A_{H_k-}(\lambda)$ ($k=1,2$) by formula (4.3) and condition (4.1):

$$A_{12+}(\lambda) = (H_1 - H_2)^{-1}\big[A_{H_2+}(\lambda) - A_{H_1+}(\lambda)\big],$$
$$A_{11-}(\lambda) = (H_1 - H_2)^{-1}\big[A_{H_1-}(\lambda)H_1 - A_{H_2-}(\lambda)H_2\big],$$
$$A_{22+}(\lambda) = A_{H_1+}(\lambda) + H_1 A_{12+}(\lambda) = A_{H_2+}(\lambda) + H_2 A_{12+}(\lambda),$$
$$A_{21-}(\lambda) = H_1 A_{11-}(\lambda) - A_{H_1-}(\lambda)H_1 = H_2 A_{11-}(\lambda) - A_{H_2-}(\lambda)H_1.$$

As already is known that the transmission matrix $\Pi(\lambda)$ is closely related with the scattering matrix for the system of ordinary differential equation on the whole axis. Indeed, if we take the coefficients zero for $x < 0$, then we obtain the system (4.5) on the whole axis and the transmission matrix for the system (1.1) is coincide with the scattering matrix on the whole axis for the system (4.5). The uniqueness of the ISP for system (4.5) with the potential (1.2) is studied in [11-14] (see also Theorem 2.2 in [15]). It means that the potential $Q(x)$ is uniquely determined from $\Pi(\lambda)$. Theorem 2 is proved.

The matrix Riemann-Hilbert problems which is mentioned in Theorem 2 are in the form of

$$[I + A_{H+}(\lambda)]S_H(\lambda) = [I + A_{H-}(\lambda)], \quad \lambda \in \mathbb{R}$$

with boundary conditions (3.4), that is,

$$A_{H\pm}(\infty) = 0$$

where $A_{H+}(\lambda)$ and $A_{H-}(\lambda)$ are $n \times n$ matrices which are analytic in upper and lower complex $\lambda$-plane and the components of $A_{H\pm}(\lambda)$ belong to set $\mathfrak{S}^{\pm}$, which is denote the set consisting of functions of the form $\int_0^{+\infty} f(x)e^{\pm i\lambda x}dx$, $\lambda \in \mathbb{R}$, where $f(x)$ is continuous and $f(x) \in L_1$, that is, $\int_0^{\pm\infty}|f(x)|\,dx$ exists.

If the matrix functions $I + A_{H+}(\lambda)$ and $I + A_{H-}(\lambda)$ degenerate nowhere in their domains of analyticity, i. e., $det[I + A_{H+}(\lambda)] \neq 0$, $det[I + A_{H-}(\lambda)] \neq 0$, then the Riemann-Hilbert problem is called to be regular. The solution of a regular Riemann-Hilbert problem under canonical normalization is unique [11, p. 155]. By matrix analogue of the Wiener theorem [25, p. 60-63], under the condition $det[I + A_{H+}(\lambda)] \neq 0$, there exists a matrix $B_{H+}(\lambda)$ with the components belonging to $\mathfrak{S}^+$, such that $[I + A_{H+}(\lambda)]^{-1} = I + B_{H+}(\lambda)$. Thus, in the regular case the matrix RH problem reduces to left canonical factorization problem of the matrix $S_H(\lambda)$ ([26], p. 31-37). Because the factorization factors are uniquely determined in the left canonical factorization problem ([26], p. 35-37), the matrices $A_{H+}(\lambda)$ and $A_{H-}(\lambda)$ are uniquely determined by $S_H(\lambda)$.

Therefore, we become the following corollary of the Theorem 2.

**COROLLARY.** *If the Riemann-Hilbert problems (4.2) are regular, then the matrix potential $Q(x)$ of the $\mathcal{M}$-canonical system (1.1) is uniquely determined from its scattering matrices $S_{H_1}(\lambda)$ and $S_{H_2}(\lambda)$ with $\det(H_1 - H_2) \neq 0$.*

**REMARK.** If the RH problem is not regular, but the matrix $S_H(\lambda)$ is nonsingular, that is, $\det S_H(\lambda) \neq 0$, the solution of considered matrix Riemann-Hilbert problem yields the solution to following system of singular integral equation by using the Plemelj formula ([27], p. 589):

$$r(\lambda)T(\lambda) = I + \frac{1}{2\pi i} \oint_{-\infty}^{+\infty} \frac{r(s)g(s)}{s - \lambda} ds$$

where $r(\lambda) = I + A_{H+}(\lambda)$, $g(\lambda) = S_H(\lambda) - I$, $T(\lambda) = I + \frac{1}{2}g(\lambda)$. It is shown in [28] that in the case where $g(\lambda), g'(\lambda) \in L_2 \cap L_\infty$, this equation is Fredholm integral equation of second type. For such Fredholm integral equation it is known that a solution exists so long as the only solution of the homogeneous problem is the zero solution. This fact also establishes the uniqueness of the solution because the difference of any two solution satisfies the homogeneous problem. It is shown in [28] that if either $Re\, S_H(\lambda)$ or $Im\, S_H(\lambda)$ is definite then there exists no nontrivial homogeneous solution.

For a detailed exposition of this extremely interesting area of modern complex analysis we refer the reader to the monographs [29,30].

## 1. EXAMPLE

Consider the following $\mathcal{M}$ - canonical system on the half-axis $x \geq 0$

$$\begin{cases} -iz_1' = \lambda \xi_1 z_1, \\ -iz_2' + c_{2,1}(x)z_1 + c_{2,2n}(x)z_{2n} = \lambda \xi_2 z_2, \\ \quad \vdots \\ -iz_{2n-1}' + c_{2n-1,1}(x)z_1 + c_{2n-1,2n}(x)z_{2n} = \lambda \xi_{2n-1} z_{2n-1}, \\ -iz_{2n}' = \lambda \xi_{2n} z_{2n}, \end{cases} \quad (E.1)$$

where $\xi_1 < \ldots < \xi_n < 0 < \xi_{n+1} < \ldots < \xi_{2n}$.

Consider the system (E.1) under the boundary condition

$$\begin{bmatrix} z_{n+1}(0) \\ \vdots \\ z_{2n}(0) \end{bmatrix} = H \begin{bmatrix} z_1(0) \\ \vdots \\ z_n(0) \end{bmatrix}, \text{ with } H = \begin{bmatrix} 0 \\ \vdots & H_1 = [h_{kj}]_{j=2}^n \\ 0 \\ 1 & 0 & \cdots & 0 \end{bmatrix} \text{ and } \det H_1 \neq 0.$$

(E.2)

It is easy to check that the system (E.1) with the asymptotic

$$z_k(x) = a_k e^{i\lambda \xi_k x} + o(1), \quad x \to +\infty,$$

$$z_{n+k}(x) = b_{n+k} e^{i\lambda \xi_{n+k} x} + o(1), \quad x \to +\infty, k = 1, \cdots n$$

has the explicit solution in the following form

$$z_1(x) = a_1 e^{i\lambda \xi_1 x}, \quad z_{2n}(x) = b_{2n} e^{i\lambda \xi_{2n} x},$$

$$z_k(x) = a_k e^{i\lambda \xi_k x} - ia_1 \int_x^{+\infty} c_{k,1}(t) e^{i\lambda(\xi_1-\xi_k)t} dt\, e^{i\lambda \xi_k x} - ib_{2n} \int_x^{+\infty} c_{k,2n}(t) e^{i\lambda(\xi_{2n}-\xi_k)t} dt\, e^{i\lambda \xi_k x},$$
$$k = 2, \cdots n,$$

$$z_{n+k}(x) = b_{n+k} e^{i\lambda \xi_{n+k} x} - ia_1 \int_x^{+\infty} c_{n+k,1}(t) e^{i\lambda(\xi_1-\xi_{n+k})t} dt\, e^{i\lambda \xi_{n+k} x}$$
$$- ib_{2n} \int_x^{+\infty} c_{n+k,2n}(t) e^{i\lambda(\xi_{2n}-\xi_{n+k})t} dt\, e^{i\lambda \xi_{n+k} x}, \qquad k = 1, \cdots n-1.$$

Taking into account the boundary conditions (E.2), i.e. the conditions

$$z_{2n}(0) = z_1(0),$$
$$z_{n+k}(0) = \sum_{j=2}^{n} h_{kj} z_j(0), \qquad k = 1, \cdots, n-1$$

we obtain the following relations between column vectors $A = (a_1, \cdots, a_n)$ and $B = (b_{n+1}, \cdots, b_{2n})$:

$$B = S_H(\lambda) H A,$$

where

$$S_H(\lambda) = I + \begin{bmatrix} 0 & \cdots & 0 & S_{1,n}(\lambda) \\ \vdots & \ddots & \vdots & \vdots \\ 0 & \cdots & 0 & S_{n-1,n}(\lambda) \\ 0 & \cdots & 0 & 0 \end{bmatrix},$$

$$S_{k,n}(\lambda) = \int_0^{+\infty} i \left[ \frac{1}{\xi_{n+k} - \xi_1} c_{n+k,1}\left(\frac{s}{\xi_{n+k}-\xi_1}\right) - \sum_{j=2}^{n} \frac{h_{kj}}{\xi_j - \xi_1} c_{j,1}\left(\frac{s}{\xi_j - \xi_1}\right) \right] e^{-i\lambda s} ds$$

$$+ \int_0^{+\infty} i \left[ \frac{1}{\xi_{2n} - \xi_{n+k}} c_{n+k,2n}\left(\frac{s}{\xi_{2n}-\xi_{n+k}}\right) \right. \hspace{3cm} \text{(E. 3)}$$

$$\left. - \sum_{j=2}^{n} \frac{h_{kj}}{\xi_{2n} - \xi_j} c_{j,2n}\left(\frac{s}{\xi_{2n} - \xi_j}\right) \right] e^{i\lambda s} ds, \quad k = 1, \cdots, n-1.$$

Now, consider the ISP for the system (E.1), i.e. the problem of recovering the coefficients $c_{k1}(x), c_{k,2n}(x) (k = 2, \cdots, 2n-1)$ of the system (1.1) from its scattering matrix $S_H(\lambda)$ on the half-axis. As is shown, the coefficients of the system (E.1) and scattering matrix $S_H(\lambda)$ are related with the relation (E.3).

Denoting

$$C_{k-}(\lambda) = \int_0^{+\infty} c_{k-}(s) e^{-i\lambda s} ds,$$

$$c_{k-}(s) = i \left[ \frac{1}{\xi_{n+k} - \xi_1} c_{n+k,1}\left(\frac{s}{\xi_{n+k}-\xi_1}\right) - \sum_{j=2}^{n} \frac{h_{kj}}{\xi_j - \xi_1} c_{j,1}\left(\frac{s}{\xi_j - \xi_1}\right) \right]$$

and

$$C_{k+}(\lambda) = \int_0^{+\infty} c_{k+}(s) e^{i\lambda s} ds,$$

$$c_{k+}(s) = i\left[\frac{1}{\xi_{2n}-\xi_{n+k}} c_{n+k,2n}\left(\frac{s}{\xi_{2n}-\xi_{n+k}}\right) - \sum_{j=2}^{n} \frac{h_{kj}}{\xi_{2n}-\xi_j} c_{j,2n}\left(\frac{s}{\xi_{2n}-\xi_j}\right)\right]$$

the formula (E.3) can be written in following form

$$S_{k,n}(\lambda) = C_{k+}(\lambda) + C_{k-}(\lambda). \tag{E.4}$$

By using (1.2), it is easy to show that

$$|c_{k\pm}(s)| \leq M e^{-\frac{\varepsilon}{\xi_{2n}-\xi_1} s},$$

where $M$ is constant. Then we can conclude that, the matrix function $C_{k+}(\lambda)$ is analytic in half-plane $Im\lambda > -\frac{\varepsilon}{\xi_{2n}-\xi_1}$, and $C_{k-}(\lambda)$ is analytic in half-plane $Im\lambda < \frac{\varepsilon}{\xi_{2n}-\xi_1}$. In addition the functions $C_{k\pm}(\lambda)$ tend to zero as $Im\lambda \to \infty$ in the domains of analyticity.

Thus, the ISP for the system (E.1) on the half-axis can be solvable by Wiener-Hopf method. Actually, it is possible to determine the functions $C_{k-}(\lambda)$ and $C_{k+}(\lambda)$ of a complex variable $\lambda$, which are analytic respectively in the half-plane $Im\lambda < \frac{\varepsilon}{\xi_{2n}-\xi_1}$ and $Im\lambda > -\frac{\varepsilon}{\xi_{2n}-\xi_1}$, tend to zero as $Im\lambda \to \infty$ in both domains of analyticity and satisfy in the strip $-\frac{\varepsilon}{\xi_{2n}-\xi_1} < Im\lambda < \frac{\varepsilon}{\xi_{2n}-\xi_1}$ the equation (П.4).

Since the function $S_{k,n}(\lambda)$ is analytic in the strip $-\frac{\varepsilon}{\xi_{2n}-\xi_1} < Im\lambda < \frac{\varepsilon}{\xi_{2n}-\xi_1}$, then the following representation is possible in the given strip

$$S_{k,n}(\lambda) = S_{k,n+}(\lambda) + S_{k,n-}(\lambda), \tag{E.5}$$

when $S_{k,n}(\lambda)$ tends uniformly to zero as $|\lambda| \to +\infty$ in this strip (see [31, p. 293]). Here the functions $S_{k,n+}(\lambda)$ and $S_{k,n-}(\lambda)$ are analytic in $Im\lambda > -\frac{\varepsilon}{\xi_{2n}-\xi_1}$ and $Im\lambda < \frac{\varepsilon}{\xi_{2n}-\xi_1}$, respectively.

From (E.4) and (E.5) we get the following formula

$$C_{k+}(\lambda) - S_{k,n+}(\lambda) = -C_{k-}(\lambda) + S_{k,n-}(\lambda) \tag{E.6}$$

in the strip $-\frac{\varepsilon}{\xi_{2n}-\xi_1} < Im\lambda < \frac{\varepsilon}{\xi_{2n}-\xi_1}$.

The left side of (E.6) is a function which is analytic in half-plane $Im\lambda > -\frac{\varepsilon}{\xi_{2n}-\xi_1}$ and the right side of (E.6) is analytic in half-plane $Im\lambda < \frac{\varepsilon}{\xi_{2n}-\xi_1}$. From the equality of these functions in the strip $-\frac{\varepsilon}{\xi_{2n}-\xi_1} < Im\lambda < \frac{\varepsilon}{\xi_{2n}-\xi_1}$ it follows that there exist a unique entire function $p(\lambda)$ coinciding, respectively, with the left and right sides of (E.6) in the domains of their analyticity. Since the function $C_{k\pm}(\lambda) - S_{k,n\pm}(\lambda)$ tends zero as $Im\lambda \to \infty$ in the domain of analyticity, then $p(\lambda) = 0$. So that $C_{k-}(\lambda) = S_{k,n-}(\lambda)$ and $C_{k+}(\lambda) = S_{k,n+}(\lambda)$. In means that the functions $S_{k,n-}(\lambda)$ and $S_{k,n+}(\lambda)$ are Laplace transformations of the functions $c_{k-}(s)$ and $c_{k+}(s)$, respectively. Then

$$c_{k+}(s) = \frac{1}{2\pi i} \int_{-\infty+ix}^{+\infty+ix} S_{k,n+}(\lambda) e^{i\lambda s} d\lambda, \qquad x > -\frac{\varepsilon}{\xi_{2n}-\xi_1}$$

$$c_{k-}(s) = \frac{1}{2\pi i} \int_{-\infty+ix}^{+\infty+ix} S_{k,n-}(\lambda) \, e^{-i\lambda s} d\lambda, \qquad x < \frac{\varepsilon}{\xi_{2n} - \xi_1}.$$

Thus, we obtain the following system of linear algebraic equations with respect to $c_{k1}(x)$, $c_{k,2n}(x)$ ($k = 2, \cdots, 2n-1$):

$$\frac{1}{\xi_{n+k} - \xi_1} c_{n+k,1}\left(\frac{s}{\xi_{n+k} - \xi_1}\right) - \sum_{j=2}^{n} \frac{h_{kj}}{\xi_j - \xi_1} c_{j,1}\left(\frac{s}{\xi_j - \xi_1}\right) = -ic_{k-}(s),$$

(E.7)

$$\frac{1}{\xi_{2n} - \xi_{n+k}} c_{n+k,2n}\left(\frac{s}{\xi_{2n} - \xi_{n+k}}\right) - \sum_{j=2}^{n} \frac{h_{kj}}{\xi_{2n} - \xi_j} c_{j,2n}\left(\frac{s}{\xi_{2n} - \xi_j}\right) = -ic_{k+}(s),$$

$$k = 1, \cdots n - 1.$$

It easy to see that, in the case $n \geq 2$, the system (E.7) has not a unique solution, therefore the uniqueness of the ISP for the system (E.1) is violated.

Consider the new ISP for the system (E.1) on the half-axis with the boundary condition

$$\begin{bmatrix} z_{n+1}(0) \\ \vdots \\ z_{2n}(0) \end{bmatrix} = \widetilde{H} \begin{bmatrix} z_1(0) \\ \vdots \\ z_n(0) \end{bmatrix}, \text{ with } \widetilde{H} = \begin{bmatrix} 0 & & & \\ \vdots & \widetilde{H}_1 = [\tilde{h}_{kj}]_{j=2}^{n} \\ 0 & & & \\ 1 & 0 & \cdots & 0 \end{bmatrix} \text{ and } \det\widetilde{H}_1 \neq 0.$$

(E.8)

It is easy to see that $\det(H - \widetilde{H}) = 0$.

Let

$$S_{\widetilde{H}}(\lambda) = I + \begin{bmatrix} 0 & \cdots & 0 & \tilde{S}_{1,n}(\lambda) \\ \vdots & \ddots & \vdots & \vdots \\ 0 & \cdots & 0 & \tilde{S}_{n-1,n}(\lambda) \\ 0 & \cdots & 0 & 0 \end{bmatrix}$$

be a scattering matrix for the system (E.1) on the half-axis with the boundary condition (E.8). Then analogically to formulas (E.7) we obtain the following formulas for the coefficients $c_{k1}(x)$, $c_{k,2n}(x)$ ($k = 2, \cdots, 2n-1$) of the system (E.1):

$$\frac{1}{\xi_{n+k} - \xi_1} c_{n+k,1}\left(\frac{s}{\xi_{n+k} - \xi_1}\right) - \sum_{j=2}^{n} \frac{\tilde{h}_{kj}}{\xi_j - \xi_1} c_{j,1}\left(\frac{s}{\xi_j - \xi_1}\right) = -i\tilde{c}_{k-}(s),$$

(E.9)

$$\frac{1}{\xi_{2n} - \xi_{n+k}} c_{n+k,2n}\left(\frac{s}{\xi_{2n} - \xi_{n+k}}\right) - \sum_{j=2}^{n} \frac{\tilde{h}_{kj}}{\xi_{2n} - \xi_j} c_{j,2n}\left(\frac{s}{\xi_{2n} - \xi_j}\right) = -i\tilde{c}_{k+}(s),$$

$$k = 1, \cdots n - 1.$$

where

$$\tilde{c}_{k+}(s) = \frac{1}{2\pi i} \int_{-\infty+ix}^{+\infty+ix} \tilde{S}_{k,n+}(\lambda) \, e^{i\lambda s} d\lambda, \qquad x > -\frac{\varepsilon}{\xi_{2n} - \xi_1},$$

$$\tilde{c}_{k-}(s) = \frac{1}{2\pi i} \int_{-\infty+ix}^{+\infty+ix} \tilde{S}_{k,n-}(\lambda) e^{-i\lambda s} d\lambda, \qquad x < \frac{\varepsilon}{\xi_{2n} - \xi_1}$$

and

$$\tilde{S}_{k,n}(\lambda) = \tilde{S}_{k,n+}(\lambda) + \tilde{S}_{k,n-}(\lambda).$$

It is easy to see that, in the case $n > 2$, the uniqueness of the solution of the system (E.7), (E.9) is violated if $det\ (H_1 - \widetilde{H}_1) = 0$. In this case, the ISP for the system (E.1) has also not a unique solution.